\documentclass[12pt]{article}
\usepackage{amsmath}
\usepackage{amsfonts}
\usepackage{latexsym}
\usepackage{amssymb}

\def\del#1{}

\title{{\scshape\normalsize Mathematics Division, National Center for
         Theoretical Sciences at Taipei}\\
         {\scshape\large NCTS/TPE-Math Technical Report 2004-011} \\~\\
         \bf Distance-two labelings of digraphs}

\author{Gerard J. Chang\thanks{Department of Mathematics,
                               National Taiwan University,
                               Taipei 10617, Taiwan.
                               Email: gjchang@math. ntu.edu.tw.
                               Supported in part by the National Science
                               Council under grant
                               NSC90-2115-M002-024.}
                       \thanks{Member of Mathematics Division, National
                               Center for Theoretical Sciences at Taipei.}
        \and
        Jer-Jeong Chen\thanks{Department of Marketing and Logistics,
                              Chung Kuo Institute of Technology,
                              Taipei 116, Taiwan.
                              Email: jungle@ms.ckitu.edu.tw.
                              Supported in part by the National Science
                              Council under grant
                              NSC-90-2115-M163-001.} $^\dag$
        \and
        David Kuo\thanks{Department of Applied Mathematics,
                         National Dong Hwa University,
                         Hualien 974, Taiwan.
                         Email: davidk@server.am.ndhu.edu.tw.
                         Supported in part by the National Science
                         Council under grant
                         NSC89-2115-M-259-007.} $^\dag$
        \and
        Sheng-Chyang Liaw\thanks{Department of Mathematics,
                                 National Central University,
                                 Chungli 32054, Taiwan.
                                 Email: scliaw@math. ncu.edu.tw.
                                 Supported in part by the National Science
                                 Council under grant
                                 NSC-90-2115-M-008-011.} $^\dag$
       }
\date{June 19, 2004}   
\date{July 8, 2004}    

\begin{document}
\maketitle
\newcommand{\w}{{\overline w}}
\newtheorem{theorem}{Theorem}
\newtheorem{lemma}[theorem]{Lemma}
\newtheorem{corollary}[theorem]{Corollary}
\newtheorem{definition}[theorem]{Definition}
\newtheorem{proposition}[theorem]{Proposition}
\def\qed{\hfill\rule{0.5em}{0.809em}\medskip}
\def\proof{\noindent{\bf Proof.}\ \ }
\baselineskip=22.1pt
\parindent=1cm
\long\def\longdelete#1{}

\begin{abstract}
For positive integers $j \ge k$,
an $L(j,k)$-labeling of a digraph $D$ is a function $f$ from $V(D)$ into
the set of nonnegative integers such that $|f(x)-f(y)|\geq j$ if $x$ is
adjacent to $y$ in $D$ and $|f(x)-f(y)|\geq k$ if $x$ is of distant two
to $y$ in $D$.
Elements of the image of $f$ are called labels.
The $L(j,k)$-labeling problem is to determine the $\vec{\lambda}_{j,k}$-number
$\vec{\lambda}_{j,k}(D)$ of a digraph $D$,
which is the minimum of the maximum label used in an $L(j,k)$-labeling of $D$.
This paper studies $\vec{\lambda}_{j,k}$-numbers of digraphs.
In particular,
we determine $\vec{\lambda}_{j,k}$-numbers of digraphs whose longest dipath
is of length at most $2$,
and $\vec{\lambda}_{j,k}$-numbers of ditrees having dipaths of length $4$.
We also give bounds for $\vec{\lambda}_{j,k}$-numbers of bipartite digraphs
whose longest dipath is of length $3$.
Finally,
we present a linear-time algorithm for determining
$\vec{\lambda}_{j,1}$-numbers of ditrees
whose longest dipath is of length $3$.

\bigskip

\end{abstract}

\newpage
\section{Introduction}

For positive integers $j \ge k$, an {\it $L(j,k)$-labeling} of a
graph $G$ is a function $f$ from $V(G)$ into the set of
nonnegative integers such that $|f(x)-f(y)|\ge j$ if $x$ is adjacent to
$y$ in $G$ and $|f(x)-f(y)|\ge k$ if $x$ is of distance two
to $y$ in $G$. Elements of the image of $f$ are called {\it
labels}, and the {\it span} of $f$ is the difference between the
largest and the smallest labels of $f$. The minimum span taken
over all $L(j,k)$-labelings of $G$, denoted by $\lambda_{j,k}(G)$,
is called the {\it $\lambda_{j,k}$-number} of $G$. And, if $f$ is
a labeling with a minimum span, then $f$ is called a {\it
$\lambda_{j,k}$-labeling} of $G$.  We shall assume without loss of
generality that the minimum label of an $L(j,k)$-labeling of $G$ is
$0$. We use $\lambda_{j}$ for $\lambda_{j,1}$ and $\lambda$ for
$\lambda_{2,1}$ for short.

A variation of Hale's channel assignment problem \cite{h1980}, the
problem of labeling a graph with a condition at distance two, was
first investigated in the case of $j=2$ and $k=1$ by Griggs and
Yeh \cite{gy1992}.  They derived formulas for the
$\lambda$-numbers of paths and cycles and established bounds on
the $\lambda$-numbers of trees and $n$-cubes.  They also
investigated the relationship between $\lambda(G)$ and other graph
invariants such as $\chi(G)$ and $\Delta(G)$. Other authors have
subsequently contributed to the literature of $L(j,k)$-labelings
with focus on the case of $j=2$ and $k=1$, see the references.

In this paper we consider the channel assignments in which
frequency inference has direction.  The formulation then becomes
$L(j,k)$-labelings on digraphs rather than on graphs. Recall that
in a digraph $D$, the {\it distance} $d_D(x,y)$ from vertex $x$ to
vertex $y$ is the length of a shortest dipath from $x$ to $y$. We
then may define $L(j,k)$-labelings for digraphs in precisely the
same way as for graphs. However, to distinguish from the notation
for graphs, we use $\vec{\lambda}_{j,k}$-number
$\vec{\lambda}_{j,k}(D)$ for a digraph $D$. We also use
$\vec{\lambda}_{j}$ for $\vec{\lambda}_{j,1}$ and $\vec{\lambda}$
for $\vec{\lambda}_{2,1}$.

The $\vec{\lambda}_{1,1}$-number of a digraph is closed related to
its oriented chromatic number.  For a digraph $D$, an {\it
oriented labeling} is a function $f$ from $V(D)$ into the set of
positive integers such that $f(x)\ne f(y)$ for $xy \in E(D)$ and
whenever an ordered pair $(p,q)$ is used for an edge $xy$ as
$(f(x), f(y))$, the ordered pair $(q,p)$ is never used for any
other edge. The {\it oriented chromatic number} $\vec{\chi} (D)$
of a digraph $D$ is the minimum size of the image of an oriented
labeling of $D$. Oriented chromatic numbers have been studied in
the literature extensively.  Notice that $\vec{\lambda}_{1,1} (D)
\le \vec{\chi}(D) -1$ for any digraph $D$.

For a tree $T$, Griggs and Yeh \cite{gy1992} showed that
$\Delta(T)+1 \le \lambda(T) \le \Delta(T)+2$; and a
polynomial-time algorithm to determine the value of $\lambda(T)$
(respectively, $\lambda_j(T)$) was given by Chang and Kuo
\cite{ck1996} (respectively, Chang et al. \cite{ckkly2000}). A
surprising result by Chang and Liaw \cite{cl2003} says that
$\vec{\lambda}(T) \le 4$ for any ditree $T$, which is an
orientation of a tree.  Suppose the largest length of a dipath in
the ditree $T$ is $\ell$.  They also proved that
$\vec{\lambda}(T)=2$ if $\ell=1$; $\vec{\lambda}(T)=3$ if
$\ell=2$; $3 \le \vec{\lambda}(T) \le 4$ if $\ell =3$; and
$\vec{\lambda}(T)=4$ if $\ell \ge 4$. Determining the exact value
of $\vec{\lambda}(T)$ for the case of $\ell=3$ left open, while
there are examples showing that $\vec{\lambda}(T)$ can be $3$ or
$4$.

The main results of this paper is to determine the exact value of
$\vec{\lambda}_{j,k}(D)$ for a digraph $D$ whose longest dipath is
of length $1$ or $2$. It is also proved that $j+k \le
\vec{\lambda}_{j,k} \le j+2k$ for a bipartite digraph whose
longest dipath is of length $3$. Finally, a linear-time algorithm
is given for determining $\vec{\lambda}_j(T)$ of a ditree $T$
whose longest dipath is of length $3$.

\section{Preliminary}

In this section, we first fix some notation and terminology, and
then derive some general propositions for $\lambda_{j,k}$-numbers
of digraphs.

For a graph $G$, if $D$ is the digraph resulting from $G$ by
replacing each edge $\{x,y\}$ by two (directed) edges $xy$ and
$yx$, then $\vec{\lambda}_{j,k}(D)=\lambda_{j,k}(G)$. However,
{\it in this paper all digraphs are assumed to be strongly
simple}, i.e.~it has no loops or multiple edges or both the edges
$xy$ and $yx$.

A digraph $D$ is {\it homomorphic} to another digraph $H$ if there
is a {\it homomorphism} from $D$ to $H$, which is a function $h:
V(D) \to V(H)$ such that $xy\in E(D)$ implies $h(x)h(y)\in E(H)$.
Define
$$
     N^+_D(v)=\{u: vu \in E(D)\}, ~ ~
     N^-_D(v)=\{u: uv \in E(D)\}, ~ ~
     N_D(v)=N^+_D(v) \cup N^-_D(v).
$$
If there is no confusion on the digraph $D$, we simply use
$N^+(v)$ for $N^+_D(v)$, $N^-(v)$ for $N^-_D(v)$ and $N(v)$ for
$N_D(v)$. We call the vertices in $N^+(v)$ the {\it out-neighbors}
of $v$, in $N^-(v)$ the {\it in-neighbors} and in $N(v)$ the {\it
neighbors}. A {\it source} is a vertex with no in-neighbors, and a
{\it sink} a vertex with no out-neighbors. A {\it leaf} of a
digraph is a vertex with exactly one neighbor.

An {\it orientation} of a graph is a digraph obtained from the
graph by assigning each edge of the graph an direction.
   The {\it underlying graph} of a digraph is the graph obtained
   from the digraph by forgetting the directions of its edges.

The {\it $n$-dipath} is the digraph $\vec{P}_n$ with $V(\vec{P}_n)
= \{v_0$, $v_1$,  $\ldots$, $v_{n-1}\}$ and $E(\vec{P}_n) =
\{v_0v_1$, $v_1v_2$, $\ldots$, $v_{n-2}v_{n-1}\}$.  The {\it
$n$-dicycle} is the digraph $\vec{C}_n$ with $V(\vec{C}_n) =
\{v_0$, $v_1$, $\ldots$, $v_{n-1}\}$  and $E(\vec{C}_n) =
\{v_0v_1$, $v_1v_2$, $\ldots$, $v_{n-2}v_{n-1}$, $v_{n-1}v_0\}$.
The $n$-{\it path} $P_n$ is the underlying graph of the $n$-dipath
$\vec{P}_n$, and the $n$-{\it cycle} $C_n$ is the underlying graph
of the $n$-dicycle $\vec{C}_n$.
A {\it ditree} is an orientation of a tree. Notice that a
nontrivial ditree has at least two leaves. A digraph is {\it
bipartite} if and only if its underlying graph is bipartite.

\begin{lemma} \label{1}
If $D$ is a subdigraph of a digraph $H$, then
$\vec{\lambda}_{j,k}(D) \le \vec{\lambda}_{j,k}(H)$.
\end{lemma}

\proof The lemma follows from the fact that the restriction of an
$L(j,k)$-labeling of $H$ on $V(D)$ is an $L(j,k)$-labeling of $D$,
since $1 \le d_H(x,y) \le d_D(x,y)$ for any two distinct vertices
$x$ and $y$ in $D$. \qed

Given $n$ digraphs $D_{1},D_{2},\cdots ,D_{n},$ the {\it union} of
these $n$ digraphs, denoted by $\bigcup_{i=1}^{n}D_{i},$ is the
digraph $D$ with $V(D)=\bigcup_{i=1}^{n}V(D_{i})$ and
$E(D)=\bigcup_{i=1}^{n}E(D_{i})$. The following lemma is obvious.

\begin{lemma} \label{new}
If $D=\bigcup_{i=1}^{n}D_{i},$ then
$\vec{\lambda}_{j,k}(D)=\max\limits_{1\leq i\leq
n}\vec{\lambda}_{j,k}(D_{i}).$
\end{lemma}

\begin{lemma} \label{2}
If digraph $D$ is homomorphic to digraph $H$, then
$\vec{\lambda}_{j,k}(D) \le \vec{\lambda}_{j,k}(H)$.
\end{lemma}

\proof The lemma follows from the fact that the composition of a
homomorphism $h$ from $D$ to $H$ with an $L(j,k)$-labeling of $H$
is an $L(j,k)$-labeling of $D$, since $H$ being strongly simple
implies that $1 \le d_H(h(x),h(y)) \le d_D(x,y)$ for any two
vertices $x$ and $y$ in $D$ with $1 \le d_D(x,y) \le 2$. \qed


\begin{lemma} \label{4}
The following statements hold for any digraph $D$.
\begin{description}
\item $(1)$   $\vec{\lambda}_{j,k}(D)=0$ if and only if $D$ has no
edge.

\item $(2)$   If $D$ has at least one edge, then
$\vec{\lambda}_{j,k}(D) \ge j$.

\item $(3)$   For any digraph $D$ with at least one edge,
        $\vec{\lambda}_{j,k}(D)=j$ if and only if every vertex is either a source
        or a sink.
\end{description}
\end{lemma}

\proof
Statements (1) and (2) are trivial.

(3) By statement (2), we have $\vec{\lambda}_{j,k}(D) \ge j$. If
every vertex of $D$ is either a source or a sink, then consider
the mapping $f$ defined by $f(x)=0$ if $x$ is a source and
$f(x)=j$ otherwise. It is easy to check that $f$ is an
$L(j,k)$-labeling of $D$, which implies that
$\vec{\lambda}_{j,k}(D) \le j$ and so $\vec{\lambda}_j(D)=j$. On
the other hand, if $D$ has a vertex $y$ which is neither a source
nor a sink, then choose $x\in N^-(y)$ and $z\in N^+(y)$. For any
$L(j,k)$-labeling $g$, we have $g(x)\ne g(z)$ and $g(y)$ differs
from $g(x)$ and $g(z)$ by at least $j$. These imply that $g$ must
use a label greater than $j$, i.e., $\vec{\lambda}_{j,k}(D) > j$.
\qed

Notice that in general
$\vec{\lambda}_{j,k}(\vec{P}_n)=\lambda_{j,k}(P_n)$ and
$\vec{\lambda}_{j,k}(\vec{C}_n)=\lambda_{j,k}(C_n)$. For the
purpose of this paper, we need the values \cite{gm1995}:
$\vec{\lambda}_{j,k}(\vec{P}_3) = \vec{\lambda}_{j,k}(\vec{P}_4) =
j+k$, $\vec{\lambda}_{j,k}(\vec{P}_5) = \min\{2j,j+2k\}$,
$\vec{\lambda}_{j,k}(\vec{C}_3) = 2j$ and
$\vec{\lambda}_{j,k}(\vec{C}_4) = j+2k$.

\section{Digraphs with a specified longest dipath length}

This section investigates digraphs in which
the length $\ell$ of a longest dipath is at most $3$.

The case when $\ell=1$ is a consequence of Lemma \ref{4} (3), as a
longest dipath of a digraph is of length 1 if and only if the
digraph has at least one edge and every vertex is either a source
or a sink.

We now consider the case when $\ell=2$.
There are two subcases.
We first deal with the case when $D$ is a bipartite digraph.

\begin{theorem} \label{12}
For any bipartite digraph $D$ whose longest dipath has length $2$,
we have $\vec{\lambda}_{j,k}(D) = j+k$
\end{theorem}

\proof According to Lemma \ref{1}, $\vec{\lambda}_{j,k}(D) \ge
\vec{\lambda}_{j,k}(\vec{P}_3)=j+k$ since $D$ has a dipath of
length $2$.

On the other hand, choose a bipartition $A\cup B$ of $V(D)$.
Define function $f$ on $V(D)$ by
 \begin{eqnarray*}
  f(x) &=& \left\{ \begin{array}{ll}
            0,   & \mbox{if } x \in A-S; \\
            k,   & \mbox{if } x \in A\cap S; \\
            j,   & \mbox{if } x \in B\cap S; \\
            j+k, & \mbox{if } x \in B-S.
                    \end{array}
            \right.
  \end{eqnarray*}
We shall check that $f$ is an $L(j,k)$-labeling of $D$ as follows,
which gives $\vec{\lambda}_{j,k} \le j+k$.

If $d_D(x,y)=1$, then $x$ and $y$ are in different parts, and $y$
is not a source.  In other words, either $x\in A$ with $y\in B-S$
or $x\in B$ with $y\in A-S$. Then, $|f(x)-f(y)| \ge j$.

If $d_D(x,y)=2$, then there is a dipath $x,w,y$.  First, $x$ and
$y$ are in the same part, and $y$ is not a source. Second, suppose
$x$ is not a source, i.e. $x$ has an in-neighbor $z$. Since $D$ is
strongly simple, $z\ne w$; and since $D$ is bipartite, $z\ne y$.
These give a dipath $z,x,w,y$ of length $3$, which is impossible.
So, $x$ is a source. Therefore, either $x\in A\cap S$ with $y\in
A-S$ or $x\in B\cap S$ with $y\in B-S$. In either case,
$|f(x)-f(y)| \ge k$.

Thus, $f$ is an $L(j,k)$-labeling of $D$ and so
$\vec{\lambda}_{j,k}(D) \le j+k$. \qed

We next consider the case when $\ell=2$ and $D$ is not bipartite.

\begin{theorem} \label{14}
For any non-bipartite digraph $D$ whose longest dipath is of
length $2$, we have $\vec{\lambda}_{j,k}(D) = 2j$
\end{theorem}

\proof According to Lemma \ref{new}, we may assume that $D$ is
connected.

For the case when $D$ is cyclic, since $D$ is strongly simple,
$D=\vec{C}_{3}$ for otherwise there is a dipath of length $3$.  In
this case, $\vec{\lambda}_{j,k}(D) =
\vec{\lambda}_{j,k}(\vec{C}_3)=2j$.

Now, suppose $D$ is acyclic.  Let $S_1$ denote the set of all
sources, $S_2$ all vertices which are neither a source nor a sink,
and $S_3$ all sinks. Then $S_1 \cup S_2 \cup S_3$ is a partition
of $V(D)$. Define function $f$ on $V(D)$ by $f(x)=(p-1)j$ for $x\in
S_p$.
As $|f(x)-f(y)| \ge j$ for $x\in S_p$ and $y\in S_q$ with $p\ne
q$, in order to check that $f$ is an $L(j,k)$-labeling, we only
need to show that any $S_p$ can not contain two distinct vertices
 $x$ and $y$ of distance one or two. Suppose to the contrary that
such $S_p$, $x$ and $y$ exist.  It is then obvious that $p=2$. By
the definition of $S_2$, we may choose an in-neighbor $w$ of $x$
and an out-neighbor $z$ of $y$.  Since $D$ is acyclic, $w P z$ is
a dipath of length at least $3$, which is impossible. Therefore,
$f$ is an $L(j,k)$-labeling of $D$ and so $\vec{\lambda}(D) \le
2j$.

On the other hand, suppose $D$ has an $L(j,k)$-labeling $f$ using
labels in $\{0,1,\ldots, 2j-1\}$.
 Choose an odd cycle (not necessarily directed) $(v_0, v_1,
v_2, \ldots, v_{n-1})$ in $D$. For any $p$, one of $f(v_p)$ and
$f(v_{pi+1})$ must be in $S=\{0,1,\ldots,j-1\}$ and the other in
$B=\{j,j+1,\ldots,2j-1\}$. But this is impossible as $n$ is odd.
Therefore, $\vec{\lambda}_{j,k}(D) \ge 2j$. \qed

For the case when $\ell=3$,
we only consider bipartite digraphs.

\begin{theorem} \label{15}
For any bipartite digraph $D$ whose longest dipath has length $3$,
we have $j+k \le \vec{\lambda}_{j,k}(D) \le j+2k$
\end{theorem}

\proof According to Lemma \ref{1}, $\vec{\lambda}_{j,k}(D) \ge
\vec{\lambda}_{j,k}(\vec{P}_4)=j+2k$ since $D$ has a dipath of
length $3$.

On the other hand, according to Lemma \ref{new}, we may assume
that $D$ is connected. If $D$ contains a $\vec{C}_4$, then
$D=\vec{C}_4$ for otherwise there is a dipath of length 4.
Therefor,
$\vec{\lambda}_{j,k}(D)=\vec{\lambda}_{j,k}(\vec{C}_4)=j+2k$. Now,
we may that $\vec{C}_4\not\subseteq D$. Denote by $S$ the set of
all sources and all vertices whose in-neighbors are all sources.
Choose a bipartition $A\cup B$ of $V(D)$. Then, define function
$f$ on $V(D)$ by
 \begin{eqnarray*}
  f(x) &=& \left\{ \begin{array}{ll}
            0,    & \mbox{if } x \in A-S; \\
            k,    & \mbox{if } x \in A\cap S; \\
            j+k,  & \mbox{if } x \in B\cap S; \\
            j+2k, & \mbox{if } x \in B-S.
                    \end{array}
            \right.
  \end{eqnarray*}
We check that $f$ is an $L(j,k)$-labeling of $D$ as follows.
If $|f(x)-f(y)| <j$ for some $x$ adjacent to $y$, then $x$ and $y$ are
in a same part, which contradicts that $D$ is bipartite. If
$|f(x)-f(y)|<k$ for some $x$ and $y$ with $d_D(x,y)=2$, then $x$ and
 $y$ are both in one of the sets $A-S$, $A\cap S$, $B\cap S$ and $B-S$.
By the definition of $S$, condition $d_D(x,y)=2$ implies $y\not\in
S$. This in turn implies that $x\not\in S$. Again, by the
definition of $S$, there is a vertex $z$ whose distance to $x$ is
$2$.  Since, $D$ is strongly simple, bipartite and contains no
$\vec{C}_4$, the vertices $z,x,y$ then creates a dipath of length
$4$, a contradiction. Thus, $f$ is an $L(j,k)$-labeling of $D$.
These prove the theorem. \qed

\section{Ditrees}

In this section, we studies $\vec{\lambda}_{j,k}$-numbers for
ditrees. According to the results in the previous section, we only
consider the case when the longest dipath is of length at least
$3$.

First, a useful lemma.

\begin{lemma} \label{9}
If $T$ is a ditree and $n\ge 3$,
then there is a homomorphism from $T$ to $\vec{C}_n$.
\end{lemma}

\proof The lemma is trivial when $T$ has exactly one vertex.
Suppose $T$ has at least two vertices. Choose a leaf $x$ with
(necessarily) exactly one neighbor $y$. By the induction
hypothesis, there is a homomorphism $h$ from $T-x$ to $\vec{C}_n$.
Suppose $h(y)=v_i$. We may extend $h$ to a homomorphism from $T$
to $\vec{C}_n$ by letting $h(x)=v_{i+1}$ if $x\in N^+(y)$ and
$h(x)=v_{i-1}$ if $x\in N^-(y)$, where the addition/substraction
in $i+1$ or $i-1$ are taken modula $n$. The lemma then follows
from induction. \qed

\begin{theorem} \label{11}
For any ditree $T$, we have $\vec{\lambda}_{j,k}(T) \le \min\{2j,
j+2k\}$. Moreover, $\vec{\lambda}_{j,k}(T) = \min\{2j, j+2k\}$ if
$T$ has a dipath of length $4$.
\end{theorem}

\proof According to Lemma \ref{9}, there is a homomorphism from
$T$ to $\vec{C}_3$ (respectively, $\vec{C}_4$). By Lemma \ref{2},
$\vec{\lambda}_{j,k}(T) \le \vec{\lambda}_{j,k}(\vec{C}_3)= 2j$
(respectively, $\vec{\lambda}_{j,k}(T) \le \vec{\lambda}_{j,k}
(\vec{C}_4)= j+2k$). On the other hand, suppose $T$ has a dipath
of length $4$. By Lemma \ref{1}, $\vec{\lambda}_{j,k}(T) \ge
\vec{\lambda}_{j,k}(\vec{P}_5)=\min\{2j,j+2k\}$ and so
$\vec{\lambda}_j(T) = \min\{2j,j+2k\}$. \qed

So far, we have determined $\vec{\lambda}_{j,k}$-numbers for all
ditrees except for the case when its longest dipath is of length
$3$. In the rest of this paper, we shall give an algorithm to
determine the value of $\vec{\lambda}_j(T)$ for a ditree $T$ whose
longest dipath is of length $3$. In this case, according to
Theorem \ref{15} either $\vec{\lambda}_j(T)=j+1$ or $j+2$. Below
are two examples showing that the two possibilities happen.
Consider the ditree $T_1=\vec{P}_4$
and the ditree $T_2=(V_2,E_2)$ with
$$
    V_2=\{v_1,v_2,v_3,v_4,v_5,v_6,v_7,v_8\} {\rm \ \ and \ \ }
    E_2=\{v_1v_2, v_2v_3, v_3v_4, v_5v_4, v_5v_6, v_6v_7, v_7v_8\}.
$$
It is the case that the longest dipaths of both $T_1$ and $T_2$
are of length $3$, but $\vec{\lambda}_j(T_1)=j+1$ while
$\vec{\lambda}_j(T_2)=j+2$ when $j\ge 2$.

In sprite of the algorithm for $\lambda$-nmbers of trees
\cite{ck1996}, we now give an algorithm to determine if a ditree
has an $L(j,1)$-labeling of span $j+1$, without any assumption on
length of a longest  dipath in it. For this purpose, we consider
$T$ as rooted at a vertex $v$ and denote the ditree as $T_v$ if
necessary. Let $T_v^+$ (respectively, $T_v^-$) denote the ditree
obtained from $T_v$ by adding a new vertex $v^+$ (respectively,
$v^-$) and a new edge $v v^+$ (respectively, $v^- v$). We consider
$T_v^+$ as rooted at $v^+$ and $T_v^-$ rooted at $v^-$. Denote
$$
\mbox{ $S(T_v^+)=\{(a,b): a=f(v^+)$ and $b=f(v)$ for some
       $L(j,1)$-labeling $f$ of span $j+1$ for $T_v^+\}$, }
$$
$$
\mbox{ $S(T_v^-)=\{(a,b): a=f(v^-)$ and $b=f(v)$ for some
       $L(j,1)$-labeling $f$ of span $j+1$ for $T_v^-\}$. }
$$
Note that $\vec{\lambda}_j (T_v^+)\le j+1$ if and only if
$S(T_v^+)\ne\emptyset$, and $\vec{\lambda}_j (T_v^-)\le j+1$ if
and only if $S(T_v^-)\ne\emptyset$. Suppose $T'$ is a ditree for
which we want to determine if $\vec{\lambda}_j (T')\le j+1$.
Choose any leaf $y$ in $T'$ who only neighbor is $v$ and let $T$
be the ditree obtained from $T'$ by deleting $y$ and the edge
incident with $y$. If we view $T$ as rooted at $v$, then $T'$ is
equal to either $T_v^+$ (when $vy$ is an edge in $T'$) or $T_v^-$
(when $yv$ is an edge in $T'$).

We note that $S(T_v^+)$ and $S(T_v^-)$ are subsets of the set
$$
   W = \{ (0,j), (0, j+1), (1, j+1), (j,0), (j+1,0), (j+1,1) \}.
$$
For any $a\in\{0,1,j,j+1\}$ and $S \subseteq W$, let $S_a=\{b:
(a,b)\in S\}$. Notice that $S_a$ is of size at most $2$. In fact
$S_a \subseteq W_a$ and
$$
    W_0=\{j,j+1\},
    W_1=\{j+1\},
    W_j=\{0\},
    W_{j+1}=\{0,1\}.
$$
For any $(a,b)\in W$ and $S \subseteq W$, let $S_{(a,b)}=\{b':
(a,b')\in S-\{(a,b)\}\}$. Notice that $S_{(a,b)}$ is of size at
most $1$. In fact $S_{(a,b)} \subseteq W_{(a,b)}$ and
$$
    W_{(0,j)}=\{j+1\},
    W_{(0,j+1)}=\{j\},
    W_{(1,j+1)}=W_{(j,0)}=\emptyset,
    W_{(j+1,0)}=\{1\},
    W_{(j+1,1)}=\{0\}.
$$

Suppose $T_v -v$ contains
$r+s$ ditrees $T_{u_1}, T_{u_2}, \ldots, T_{u_r},
               T_{v_1}, T_{v_2}, \ldots, T_{v_s}$,
where each $u_i$ is {\it adjacent to} $v$
and   each $v_i$ is {\it adjacent from} $v$ in $T_v$.
Note that $T_v$ can be considered as identifying
$u_1^+, u_2^+, \ldots, u_r^+, v_1^-, v_2^-, \ldots, v_s^-$ to a vertex $v$
on the disjoint union of  $T_{u_1}^+, T_{u_2}^+, \ldots, T_{u_r}^+,
                           T_{v_1}^-, T_{v_2}^-, \ldots, T_{v_s}^-$.
We then have

\begin{theorem}
For any ditree $T$, we have
$S(T_v^+) =$
$$
    \{(a,b)\in W: S(T_{u_p}^+)_{(b,a)} \ne \emptyset
                  \mbox{ for }  1\le p\le r;
                  \emptyset \ne S(T_{v_q}^-)_b
                  \ne           S(T_{u_1}^+)_{(b,a)}
                  \mbox{ for }  1\le q \le s\},
$$
where $S(T_{u_1}^+)_{(b,a)}$ is assume to be $\emptyset$ if $r=0$.
Also,
$S(T_v^-) =$
$$
    \{(a,b)\in W: S(T_{v_q}^-)_{(b,a)} \ne \emptyset
                             \mbox{ for }  1\le q\le s;
                             \emptyset \ne S(T_{u_i}^-)_b
                             \ne           (T_{v_1}^+)_{(b,a)}
                             \mbox{ for }  1\le p \le r\},
$$
where $S(T_{v_1}^+)_{(b,a)}$ is assume to be $\emptyset$ if $s=0$.
\end{theorem}

\proof The first equality follows from that fact that $T_v^+$ has
an $L(j,1)$-labeling $f$ with $f(v^+)=a$ and $f(v)=b$ if and only
if $T_{u_p}^+$ has an $L(j,1)$-labeling $g_p$ with $g_p(u_p^+)=b$
but $g_i(u_p)\ne a$ for $1\le p\le r$, and $T_{v_q}^-$ has an
$L(j,1)$-labeling $h_q$ with $h_q(v_q^-)=b$ but $h_q(v_q)\ne
g_1(u_1)$ for $1\le q\le s$.

The second equality follows from a similar argument.
\qed

We remark that in the calculation of $S(T_v^+)$, the condition
``$S(T_{v_q}^-)_b \ne S(T_{u_1}^+)_{(b,a)}$
                for $1\le q\le s$''
is redundant if $r=0$. When $r\ge 1$ and
$S(T_{u_p}^+)_{(b,a)} \ne \emptyset$  for $1\le q\le s$,
it is the case that these sets are equal to a same set.
Similar statements are true for the set $S(T_v^-)$.

Having the theorem above, we then can compute the sets $S(T_v^+)$
and $S(T_v^-)$ recursively, using the initial conditions that
$S(T_v^+)=S(T_v^-)=W$ when $T$ is a ditree of just one vertex.
This gives a linear-time algorithm to determine whether
$\vec{\lambda}_j(T)=j+1$ or not for a ditree $T$, without any
assumption on the length of the longest dipath in $T$.

\bigskip

\noindent {\bf Acknowledgements.} The authors thank Andre Rapand
for bringing attention the relation between
$\vec{\lambda}_{j,k}$-numbers and oriented chromatic numbers, and
suggesting a better notation $\vec{\lambda}_{j,k}$ while our
original notation was $\lambda^*_{j,k}$.

\def\del#1\par{}
\def\has{}
\def\manu{}
\def\new{}

\end{document}